\documentclass[letterpaper, 10 pt, conference]{ieeeconf}

\IEEEoverridecommandlockouts

\title{\LARGE \bf
	Numerical optimal control for distributed delay differential equations:\\A simultaneous approach based on linearization of the delayed variables
}

\author{Tobias K. S. Ritschel
	\thanks{*T. K. S. Ritschel is with Department of Applied Mathematics and Computer Science, Technical University of Denmark, DK-2800 Kgs. Lyngby, Denmark
		{\tt\small tobk@dtu.dk}}%
}

\usepackage{amsmath}
\usepackage{amssymb}
\usepackage{xcolor}
\usepackage{graphicx}
\usepackage{epstopdf}
\usepackage{tikz}
\usepackage{cite}

\usetikzlibrary{positioning}
\usetikzlibrary{calc}
\usetikzlibrary{shapes.geometric}
\usetikzlibrary{angles}
\usetikzlibrary{decorations.pathreplacing}
\usetikzlibrary{calligraphy}

\newcommand{\diff}[2]{\frac{\incr{#1}}{\incr{#2}}}
\newcommand{\pdiff}[2]{\frac{\partial{#1}}{\partial{#2}}}

\newcommand{\incr}{\,\mathrm{d}}

\begin{document}
	\maketitle
\thispagestyle{empty}
\pagestyle{empty}

	\begin{abstract}
	Time delays are ubiquitous in industrial processes, and they must be accounted for when designing control algorithms because they have a significant effect on the process dynamics. Therefore, in this work, we propose a simultaneous approach for numerical optimal control of delay differential equations with \emph{distributed} time delays. Specifically, we linearize the delayed variables around the current time, and we discretize the resulting implicit differential equations using Euler's implicit method. Furthermore, we transcribe the infinite-dimensional optimal control problem into a finite-dimensional nonlinear program, which we solve using Matlab's \texttt{fmincon}. Finally, we demonstrate the efficacy of the approach using a numerical example involving a molten salt nuclear fission reactor.
\end{abstract}

	\section{Introduction}\label{sec:introduction}
Many industrial processes involve time delays~\cite{Kolmanovskii:Myshkis:1999}. They can arise due to flow in a pipe (advection), slow mixing (diffusion), or communication delays, and they have a significant impact on the transient behavior~\cite{Niculescu:Gu:2004}. Consequently, neglecting the time delays can lead to poor performance of model-based algorithms, e.g., for process monitoring, control, and optimization. Typically, such algorithms are based on models that consist of differential equations, and it is common to model the time delays as \emph{absolute}. In that case, the change in the current state depends on states at discrete times in the past. For instance, if the time delay is caused by flow in a pipe, a uniform velocity profile (i.e., plug flow) leads to an absolute time delay. In contrast, for Hagen-Poiseuille flow, where the velocity profile is nonuniform~\cite{Kessler:Greenkorn:1999, Bird:etal:2002}, the time delay is \emph{distributed}~\cite{Smith:2011, Hale:Lunel:1993, Kolmanovskii:Myshkis:1992}, i.e., the change in the current state depends on all past states in a given time interval. Specifically, the dependency is described by a convolution of the past states and a \emph{kernel} (also called a memory function), which assigns a relative weight to each time in the past.

Delay differential equations (DDEs) with distributed time delays have been used to model many different processes~\cite{Kolmanovskii:Myshkis:1999}, e.g., pharmacokinetics and pharmacodynamics (PK/PD)~\cite{Hu:etal:2018}, chemotherapy-induced myelosuppression (reduced bone marrow activity)~\cite{Krzyzanski:etal:2018}, the Mackey-Glass system~\cite{Nevermann:Gros:2023, Zhang:Xiao:2016}, which has applications in both respiratory and hematopoietic diseases, population models~\cite{Cassidy:etal:2019}, and pollution in fisheries~\cite{Bergland:etal:2023}. Mechanical~\cite{Aleksandrov:etal:2023} and economic~\cite{Guerrini:etal:2020} models have also been proposed.
Additionally, Rahman et al.~\cite{Rahman:etal:2015} studied the stability of networked systems with distributed delays, Yuan and Belair~\cite{Yuan:Belair:2011} present general stability and bifurcation results, Bazighifan et al.~\cite{Bazighifan:etal:2019} analyze oscillations in higher-order DDEs with distributed delays, and Cassidy~\cite{Cassidy:2021} demonstrate the equivalence between cyclic ODEs and a scalar DDE with a distributed time delay.
Despite the many applications, most theory on optimal control of DDEs with distributed time delays is developed for linear systems. For instance, Kushner and Barnea~\cite{Kushner:Barnea:1970} use policy iteration to derive an optimal control law for a finite-horizon optimal control problem (OCP), Santos et al.~\cite{Santos:etal:2009} propose an iterative approach to approximate the optimal control law of an infinite-horizon OCP, and in a later paper, Ortega-Mart{\'{i}}nez et al.~\cite{Ortega:Martinez:etal:2023} present the solution to the same infinite-horizon OCP. Furthermore, the book by O{\u{g}}uzt{\"{o}}reli~\cite{Oguztoreli:1966} presents theoretical results related to optimal control of both linear and nonlinear DDEs.
However, to the best of the author's knowledge, numerical optimal control of nonlinear DDEs with distributed time delays has not previously been considered.

In this work, we present a simultaneous approach for approximating the solution to finite-horizon OCPs involving nonlinear DDEs with distributed time delays. We approximate the DDEs by linearizing the delayed variables (e.g., the states) around the current time. The resulting approximate system is a set of implicit differential equations, and we present the stability criteria for both the original DDEs and the approximate system. Next, we discretize the implicit differential equations using Euler's implicit method, and we transcribe the infinite-dimensional OCP to a finite-dimensional nonlinear program (NLP). We use Matlab's \texttt{fmincon} to approximate the solution to the NLP. Furthermore, we derive the kernel for distributed time delays arising from flow in a pipe with continuously differentiable velocity profiles that 1)~only depend on the radial coordinate and 2)~are strictly monotonically decreasing towards the wall of the pipe. Finally, we demonstrate the efficacy of the proposed simultaneous approach using a numerical example involving a molten salt nuclear fission reactor where the salt is circulated through a heat exchanger outside of the core. We assume Hagen-Poiseuille flow in the external pipe, and we approximate the solution to an optimal power ramping problem (i.e., a tracking problem with a time-varying setpoint).

The remainder of the paper is structured as follows. In Section~\ref{sec:optimal:control:problem}, we present the OCP, and in Section~\ref{sec:linearization}, we describe the linearization of the delayed state. We present the simultaneous approach in Section~\ref{sec:simultaneous:approach}, and in Section~\ref{sec:vel}, we derive the kernel corresponding to a nonuniform velocity profile. In Section~\ref{sec:nuclear:fission}, we present the model of the molten salt reactor, and the numerical example is presented in Section~\ref{sec:numerical:example}. Finally, we present conclusions in Section~\ref{sec:conclusions}.
	\section{Optimal control problem}\label{sec:optimal:control:problem}
We consider OCPs in the form
\begin{subequations}\label{eq:ocp}
	\begin{align}\label{eq:ocp:obj}
		\min_{\{u_k\}_{k=0}^{N-1}} \quad \phi &= \phi_x + \phi_{\Delta u},
	\end{align}
	subject to
	\begin{align}
		\label{eq:ocp:ic}
		x(t) &= x_0(t),                            \quad t \in (-\infty, t_0], \\
		\label{eq:ocp:dde}
		\dot x(t) &= f(x(t), z(t), u(t), d(t), p), \quad t \in [t_0, t_f], \\
		\label{eq:ocp:z}
		z(t) &=
		\begin{bmatrix}
			(\alpha_1 * r_1)(t) \\
			\vdots \\
			(\alpha_m * r_m)(t)
		\end{bmatrix}, \quad t \in [t_0, t_f], \\
		\label{eq:ocp:r}
		r_i(t) &= h_i(x(t), p), \quad i = 1, \ldots, m, \quad t \in [t_0, t_f], \\
		\label{eq:ocp:u}
		u(t) &= u_k,      \quad t \in [t_k, t_{k+1}[, \quad k = 0, \ldots, N-1, \\
		\label{eq:ocp:d}
		d(t) &= d_k, \quad t \in [t_k, t_{k+1}[, \quad k = 0, \ldots, N-1, \\
		\label{eq:ocp:bounds:x}
		x_{\min} &\leq x(t) \leq x_{\max}, \quad t \in [t_0, t_f], \\
		\label{eq:ocp:bounds:u}
		u_{\min} &\leq u_k \leq u_{\max}, \quad k = 0, \ldots, N-1,
	\end{align}
\end{subequations}
where the objective function in~\eqref{eq:ocp:obj} consists of a Lagrange term and an input rate-of-movement penalization term:
\begin{subequations}
	\begin{align}
		\label{eq:ocp:obj:x}
		\phi_x &= \int_{t_0}^{t_f} \Phi(x(t), u(t), d(t), p) \incr t, \\
		\label{eq:ocp:obj:u}
		\phi_{\Delta u} &= \frac{1}{2} \sum_{k=0}^{N-1} \frac{\Delta u_k^T W_k \Delta u_k}{\Delta t}.
	\end{align}
\end{subequations}
Here, $t$ is time, $t_0$ and $t_f$ are the initial and final time of the prediction and control horizon, $x$ is the state, $x_0$ is the initial state function, $z$ is the memory state given by~\eqref{eq:ocp:z}, $\alpha_i$ is the $i$'th kernel, $r_i$ is the $i$'th delayed variable given by~\eqref{eq:ocp:r}, $m$ is the number of time delays, $u$ is the manipulated input, $d$ is the disturbance variable, and $p$ are the parameters. Furthermore, \eqref{eq:ocp:ic} is an initial condition, \eqref{eq:ocp:dde} is the DDE, \eqref{eq:ocp:u}--\eqref{eq:ocp:d} are zero-order-hold (ZOH) parametrizations of the manipulated inputs and disturbance variables, and~\eqref{eq:ocp:bounds:x}--\eqref{eq:ocp:bounds:u} are bound constraints on the states and the manipulated inputs. The decision variables are the piecewise constant manipulated inputs in the $N$ control intervals. Furthermore, the kernels integrate to one, they depend on the manipulated inputs, and the convolutions in~\eqref{eq:ocp:z} are given by
\begin{align}\label{eq:convolution}
	(\alpha_i * r_i)(t) &= \int_{-\infty}^t \alpha_i(t - s, u(t)) r_i(s) \incr s.
\end{align}
Finally, $\Phi$ in~\eqref{eq:ocp:obj:x} is the stage cost, and the change in the manipulated inputs in~\eqref{eq:ocp:obj:u} is
\begin{align}
	\Delta u_k &= u_k - u_{k-1},
\end{align}
where the reference input $u_{-1}$ is given, $W_k$ is a symmetric positive definite weight matrix, and $\Delta t = t_{k+1} - t_k$ is the size of the control intervals.

\subsection{Steady state and stability}\label{sec:optimal:control:problem:steady:state:stability}
The steady state, $x_s$, $z_s$, $u_s$, and $d_s$, of the system~\eqref{eq:ocp:dde}--\eqref{eq:ocp:r} satisfies
\begin{subequations}\label{eq:dde:steady:state}
	\begin{align}
		0 &= f(x_s, z_s, u_s, d_s, p), \\
		z_s &=
		\begin{bmatrix}
			r_{1, s} \\
			\vdots \\
			r_{m, s}
		\end{bmatrix}, \\
		r_{i, s} &= h_i(x_s, p), \quad i = 1, \ldots, m.
	\end{align}
\end{subequations}
In order to describe the stability criteria, we introduce the auxiliary variable $v_i(t) = (\alpha_i * r_i)(t)$.
The steady state is asymptotically stable if the real parts of all the roots of the characteristic function are negative~\cite{Cushing:1977}. The characteristic equation is
\begin{multline}
	\det\left(\lambda I - \pdiff{f}{x} - \sum_{i=1}^m \pdiff{f}{z} \pdiff{z}{v_i} \pdiff{r_i}{x} \int_0^{\infty} e^{-\lambda s} \alpha_i(s, u_s) \incr s\right) \\
	= 0,
\end{multline}
where $\lambda$ is the root, and the Jacobian matrices are evaluated in the steady state.
	\section{Delay linearization}\label{sec:linearization}
First, we linearize the delayed quantities around the current time:
\begin{align}\label{eq:delay:approximation}
	r_i(s) &\approx r_i(t) + \dot r_i(t) (s - t).
\end{align}
Next, we use the linearization to approximate the convolutions of the kernels and the delayed variables:
\begin{align}
	(\alpha_i * r_i)(t)
	&\approx r_i(t) \int_{-\infty}^t \alpha_i(t - s, u(t)) \incr s \nonumber \\
	&- \dot r_i(t) \int_{-\infty}^t \alpha_i(t - s, u(t)) (t - s) \incr s \nonumber \\
	&= r_i(t) - \dot r_i(t) \int_0^\infty \alpha_i(\tau, u(t)) \tau \incr \tau.
\end{align}
For the first term, we have exploited that $\alpha_i$ integrates to one, and for the second term, we have changed the integration variable to $\tau = t - s$ and switched the integration limits. Furthermore, we introduce the variable
\begin{align}\label{eq:mean}
	\gamma_i(u(t)) &= \int_0^\infty \alpha_i(\tau, u(t)) \tau \incr \tau,
\end{align}
which is analogous to the expected value of a random variable with probability density function $\alpha_i$. For some kernels, it is possible to derive an explicit expression for $\gamma_i$.
Finally, the resulting approximate system is
\begin{subequations}\label{eq:ide}
	\begin{align}
		\label{eq:ide:x}
		\dot x(t) &= f(x(t), z(t), u(t), d(t), p), \\
		\label{eq:ide:z}
		z(t) &=
		\begin{bmatrix}
			r_1(t) - \dot r_1(t) \gamma_1(u(t)) \\
			\vdots \\
			r_m(t) - \dot r_m(t) \gamma_m(u(t))
		\end{bmatrix}, \\
		\label{eq:ide:r}
		r_i(t) &= h_i(x(t), p), \quad i = 1, \ldots, m,
	\end{align}
\end{subequations}
which is a set of implicit differential equations because the time derivatives of the delayed variables appear in the expression for the memory states in~\eqref{eq:ide:z}.

\subsection{Steady state and stability}
The steady state of the approximate system~\eqref{eq:ide} is identical to the steady state of the DDEs~\eqref{eq:ocp:dde}--\eqref{eq:ocp:r}, i.e., it also satisfies~\eqref{eq:dde:steady:state}. The steady state is asymptotically stable if all roots of the characteristic function have negative real part~\cite[Prop.~2.1]{Du:etal:2013}. Each root, $\lambda$, satisfies the characteristic equation
\begin{align}
	\det\left(\lambda I - \pdiff{f}{x} - \sum_{i=1}^m \pdiff{f}{z} \pdiff{z}{v_i} \pdiff{r_i}{x} (1 - \lambda \gamma_i(u_s))\right) = 0,
\end{align}
where the auxiliary variable $v_i$ was introduced in Section~\ref{sec:optimal:control:problem:steady:state:stability}, and the Jacobian matrices are evaluated in the steady state.
	\section{Simultaneous approach}\label{sec:simultaneous:approach}
We discretize the differential equations in the approximate system~\eqref{eq:ide} using Euler's implicit method with $M$ time steps per control interval. Consequently, for each time step, the residual equations
\begin{align}\label{eq:residual:function}
	R_{k, n}
	&= R_{k, n}(x_{k, n+1}, x_{k, n}, u_k, d_k, p) \nonumber \\
	&= x_{k, n+1} - x_{k, n} - f(x_{k, n+1}, z_{k, n+1}, u_k, d_k, p) \Delta t_{k, n} \nonumber \\
	&= 0
\end{align}
must be satisfied. Furthermore, we use a backward difference approximation of the time derivatives of the delayed variables. Consequently,
\begin{align}
	\label{eq:discretized:variables:z}
	z_{k, n+1} &=
	\begin{bmatrix}
		v_{1, k, n+1} \\
		\vdots \\
		v_{m, k, n+1}
	\end{bmatrix},
\end{align}
where the auxiliary variables, $v_{i, k, n+1}$, are
\begin{align}
	\label{eq:discretized:variables:v}
	v_{i, k, n+1} &= r_{i, k, n+1} - \frac{r_{i, k, n+1} - r_{i, k, n}}{\Delta t_{k, n}} \gamma_i(u_k),
\end{align}
and the delayed variables are
\begin{align}
	\label{eq:discretized:variables:r}
	r_{i, k, n} &= h_i(x_{k, n}, p).
\end{align}
The states must be continuous across the boundaries of the control intervals and satisfy the initial condition, i.e.,
\begin{subequations}\label{eq:continuity:constraints}
	\begin{align}
		x_{0, 0} &= x_0(t_0), \\
		x_{k, 0} &= x_{k-1, M}, & k &= 1, \ldots, N-1.
	\end{align}
\end{subequations}
We use these constraints to eliminate $x_{k, 0}$ for $k = 0, \ldots, N-1$. Finally, we use a right rectangle rule to approximate the integral in the Lagrange term~\eqref{eq:ocp:obj:x} in the objective function:
\begin{align}\label{eq:obj:discretized}
	\psi_x &= \sum_{k=0}^{N-1} \sum_{n=0}^{M-1} \Phi(x_{k, n+1}, u_k, d_k, p) \Delta t_{k, n}.
\end{align}
Here, $\Delta t_{k, n} = t_{k, n+1} - t_{k, n}$.
In summary, we transcribe the OCP~\eqref{eq:ocp} into the NLP
\begin{subequations}\label{eq:nlp}
	\begin{align}\label{eq:nlp:obj}
		\min_{\{\{x_{k, n+1}\}_{n=0}^{M-1}, u_k\}_{k=0}^{N-1}} \quad & \psi = \psi_x + \phi_{\Delta u},
	\end{align}
	subject to
	\begin{align}
		\label{eq:nlp:residual:equations}
		&R_{k, n}(x_{k, n+1}, x_{k, n}, u_k, d_k, p) = 0, \\
		\label{eq:nlp:bounds:x}
		&x_{\min} \leq x_{k, n+1} \leq x_{\max}, \\
		\label{eq:nlp:bounds:u}
		&u_{\min} \leq u_k \leq u_{\max},
	\end{align}
\end{subequations}
where $k = 0, \ldots, N-1$, $n = 0, \ldots, M-1$, and the bound constraints on the states in~\eqref{eq:ocp:bounds:x} are enforced pointwise.

\subsection{Jacobian of residual functions}
For brevity of the presentation, we omit the arguments of the Jacobian matrices.
The Jacobians of the residual functions are
\begin{subequations}
	\begin{align}
		\pdiff{R_{k, n}}{x_{k, n+1}} &=	I - \left(\pdiff{f}{x} + \pdiff{f}{z} \pdiff{z_{k, n+1}}{x_{k, n+1}}\right) \Delta t_{k, n}, \\
		\pdiff{R_{k, n}}{x_{k, n}} &= -I - \pdiff{f}{z} \pdiff{z_{k, n+1}}{x_{k, n}} \Delta t_{k, n}, \\
		\pdiff{R_{k, n}}{u_k} &= -\left(\pdiff{f}{u} + \pdiff{f}{z} \pdiff{z_{k, n+1}}{u_k}\right) \Delta t_{k, n},
	\end{align}
\end{subequations}
where the Jacobian of the memory states is
\begin{align}
	\pdiff{z_{k, n+1}}{w} &=
	\begin{bmatrix}
		\pdiff{v_{1, k, n+1}}{w} \\
		\vdots \\
		\pdiff{v_{m, k, n+1}}{w}
	\end{bmatrix},
\end{align}
and $w$ represents either $x_{k, n+1}$, $x_{k, n}$, or $u_k$. Furthermore,
\begin{subequations}
	\begin{align}
		\pdiff{v_{i, k, n+1}}{x_{k, n+1}} &= \left(1 - \frac{\gamma_i(u_k)}{\Delta t_{k, n}}\right) \pdiff{r_{i, k, n+1}}{x_{k, n+1}}, \\
		\pdiff{v_{i, k, n+1}}{x_{k, n}}   &= \frac{\gamma_i(u_k)}{\Delta t_{k, n}} \pdiff{r_{i, k, n}}{x_{k, n}}, \\
		\pdiff{v_{i, k, n+1}}{u_k}        &= -\frac{r_{i, k, n+1} - r_{i, k, n}}{\Delta t_{k, n}} \pdiff{\gamma_i}{u}.
	\end{align}
\end{subequations}
Finally, the Jacobians of the delayed variables are
\begin{align}
	\pdiff{r_{i, k, n}}{x_{k, n}} &= \pdiff{h_i}{x}.
\end{align}

\subsection{Gradient of the objective function}
The gradients of the objective function are
\begin{subequations}
	\begin{align}
		\nabla_{x_{k, n+1}} \psi &= \nabla_{x_{k, n+1}} \psi_x, \\
		\nabla_{u_k} \psi &= \nabla_{u_k} \psi_x + \nabla_{u_k} \phi_{\Delta u},
	\end{align}
\end{subequations}
where the gradients of the approximate Lagrange term are
\begin{subequations}
	\begin{align}
		\nabla_{x_{k, n+1}} \psi_x &= \nabla_x \Phi (x_{k, n+1}, u_k, d_k, p) \Delta t_{k, n}, \\
		\nabla_{u_k} \psi_x &= \sum_{n=0}^{M-1} \nabla_u \Phi(x_{k, n+1}, u_k, d_k, p) \Delta t_{k, n},
	\end{align}
\end{subequations}
and the gradients of the rate-of-movement penalization term are
\begin{subequations}
	\begin{align}
		\nabla_{u_k} \phi_{\Delta u} &= \frac{W_k \Delta u_k - W_{k+1} \Delta u_{k+1}}{\Delta t}, \\
		\nabla_{u_{N-1}} \phi_{\Delta u} &= \frac{W_{N-1} \Delta u_{N-1}}{\Delta t},
	\end{align}
\end{subequations}
for $k = 0, \ldots, N-2$.

\subsection{Implementation}
We approximate the solution to the NLP~\eqref{eq:nlp} using the interior point algorithm implemented in Matlab's \texttt{fmincon}, and we supply the gradient of the objective function, and the Jacobian associated with the nonlinear equality constraints.
The Jacobian is implemented as a sparse matrix, and \texttt{fmincon} uses a finite-difference approximation of the Hessian matrix.
	\section{Distributed time delays}\label{sec:vel}
We consider a reactor with an outlet stream flowing through a cylindrical pipe of length $L$ and radius $R$. The velocity profile, $v$, only depends on the radius, $r$, it is continuously differentiable and strictly monotonically decreasing, and it is zero at the pipe wall, i.e., the flow satisfies a no-slip boundary condition such that $v(R) = 0$. The molar flow rate of component $i$ out of the reactor is
\begin{align}\label{eq:f:out}
	f_{i, out}
	&= 2 \pi \int_0^R C_i(t) v(r) r \incr r
	 = C_i(t) F,
\end{align}
where $C_i$ is the molar concentration of the $i$'th component, and $F$ is the volumetric flow rate given by
\begin{align}\label{eq:F}
	F &= 2 \pi \int_0^R v(r) r \incr r.
\end{align}
Next, we derive an expression for the molar flow rate into a subsequent reactor at the end of the pipe:
\begin{align}\label{eq:f:in}
	f_{i, in}(t)
	&= 2\pi \int_0^R C_i(t - \tau) v(r) r \incr r.
\end{align}
Here, $\tau$ is the time it takes the fluid to travel through the pipe at a given radius. We will change the integration variable in~\eqref{eq:f:in} to $\tau$, and we let the radius be an implicit function of $\tau$ defined by
\begin{align}\label{eq:velocity}
	v(r) \tau &= L.
\end{align}
The differential of the radius is
\begin{align}\label{eq:differential}
	\incr r &= \diff{r}{\tau} \incr \tau,
\end{align}
where the derivative is obtained by differentiating~\eqref{eq:velocity}:
\begin{align}
	\diff{v}{r} \diff{r}{\tau} \tau + v(r) &= 0.
\end{align}
Consequently,
\begin{align}\label{eq:radius:derivative}
	\diff{r}{\tau} &= -\frac{v(r)}{\tau} \left(\diff{v}{r}\right)^{-1} = -\frac{L}{\tau^2} \left(\diff{v}{r}\right)^{-1},
\end{align}
where we have used that $v(r) = L/\tau$. Next, we use~\eqref{eq:differential} and~\eqref{eq:radius:derivative} to change the variable of integration in~\eqref{eq:f:in}:
\begin{align}
	f_{i, in}(t) &= \int_{\tau_0}^\infty \bar \alpha(\tau) C_i(t - \tau) \incr \tau.
\end{align}
The minimum travel time is $\tau_0 = L/v(0)$, and the unnormalized kernel is
\begin{align}
	\bar \alpha(\tau) &= -2 \pi \frac{L^2}{\tau^3} \left(\diff{v}{r}\right)^{-1} r(\tau)
\end{align}
for $\tau \in [\tau_0, \infty)$ and zero otherwise. Furthermore, the integral of the unnormalized kernel is the volumetric flow rate:
\begin{align}
	\int_{\tau_0}^\infty \bar \alpha(\tau) \incr \tau &= F.
\end{align}
Therefore, the normalized kernel is
\begin{align}
	\alpha(\tau) &= \frac{\bar \alpha(\tau)}{F}.
\end{align}
Finally, we change the variable of integration to $s = t - \tau$ and use that $\mathrm d s = -\incr \tau$ and $s \rightarrow -\infty$ as $\tau \rightarrow \infty$ to obtain
\begin{align}
	f_{i, in}(t) &= F \int_{-\infty}^t \alpha(t - s) C_i(s) \incr s,
\end{align}
where we have switched the integration limits. The integral is a convolution in the form~\eqref{eq:convolution}.


\subsection{Hagen-Poiseuille flow}
For Hagen-Poiseuille flow~\cite[Sec.~6.4.2]{Kessler:Greenkorn:1999},~\cite{Bird:etal:2002}, the velocity profile is
\begin{align}\label{eq:parabolic:velocity:profile}
	v(r)     &= a (R^2 - r^2), &
	a      	 &= \frac{\Delta P}{4 \mu L},
\end{align}
where $\Delta P$ is the difference between the pipe inlet and outlet pressure, and $\mu$ is the viscosity. For this velocity profile,
\begin{align}
	\tau_0 &= \frac{L}{a R^2}, &
	F &= \frac{\pi}{2} a R^4,
\end{align}
where the expression for $F$ is the Hagen-Poiseuille equation, and
\begin{align}
	\diff{v}{r} &= -2 a r, &
	\left(\diff{v}{r}\right)^{-1} r &= -\frac{1}{2 a}.
\end{align}
Consequently, the unnormalized and normalized kernels are
\begin{align}\label{eq:hagen:poiseuille:kernel}
	\bar \alpha(\tau) &= \frac{\pi L^2}{a} \frac{1}{\tau^3}, &
	\alpha(\tau) &= 2 \frac{L^2}{a^2 R^4} \frac{1}{\tau^3} = 2 \frac{\tau_0^2}{\tau^3},
\end{align}
for $\tau \in [\tau_0, \infty)$ and zero otherwise.
Finally, the variable $\gamma$ defined in~\eqref{eq:mean} is
\begin{align}\label{eq:hagen:poiseuille:average:travel:time}
	\gamma &= \int_0^\infty \alpha(\tau) \tau \incr \tau = \int_{\tau_0}^\infty \alpha(\tau) \tau \incr \tau = 2 \tau_0.
\end{align}
	\section{Molten salt nuclear reactor}\label{sec:nuclear:fission}
In this section, we present a model of a molten salt nuclear fission reactor where the molten salt is circulated through a heat exchanger before it enters back into the reactor. The concentrations of $N_g = 6$ neutron precursor groups and the neutrons in the reactor core are described by
\begin{subequations}\label{eq:nuclear:fission:concentrations}
	\begin{align}
		\dot C_i(t) &= (C_{i, in}(t) - C_i(t)) D(t) + R_i(t), & i &= 1, \ldots, N_g, \\
		\dot C_n(t) &= R_n(t),
	\end{align}
\end{subequations}
respectively, where $n = N_g+1$, and the dilution rate,
\begin{align}
	D(t) &= F(t)/V,
\end{align}
is the ratio between the volumetric flow rate, $F$, and the volume of the reactor, $V$. The inlet concentrations of the neutron precursor groups are
\begin{align}\label{eq:nuclear:fission:inlet:flow:rate}
	C_{i, in}(t) &= e^{-\lambda_i \gamma} \int_{-\infty}^t \alpha_f(t - s) C_i(s) \incr s,
\end{align}
where the kernel $\alpha_f$ is given by the expression for $\alpha$ in~\eqref{eq:hagen:poiseuille:kernel} for the full length, $L$, of the external circulation loop and the pressure difference $\Delta P$. Furthermore, $\lambda_i$ is the decay rate of the $i$'th neutron precursor group, and for simplicity, we assume that the decay depends on the \emph{average} travel time, $\gamma$ in~\eqref{eq:hagen:poiseuille:average:travel:time}.
The production term is
\begin{align}
	R(t) &= S^T(t) r(t),
\end{align}
where $S$ is a stochiometric matrix, and $r$ is a vector of reaction rates:
{\small
\renewcommand{\arraystretch}{1.2}
\begin{align}
	S(t) &=
	\begin{bmatrix}
		-1 & & & 1 \\
		& \ddots & & \vdots \\
		& & -1 & 1 \\
		\beta_1 & \cdots & \beta_{N_g} & \rho(t) - \beta
	\end{bmatrix}\hspace{-2pt}, \;
	r(t) =
	\begin{bmatrix}
		\lambda_1 C_1(t) \\
		\vdots \\
		\lambda_{N_g} C_{N_g}(t) \\
		C_n(t)/\Lambda
	\end{bmatrix}.
\end{align}
}%
Here, $\beta_i$ is the delayed neutron fraction of precursor group $i$, and $\beta$ is the sum of $\beta_i$ for $i = 1, \ldots, N_g$. Furthermore, $\Lambda$ is the mean neutron generation time, and the reactivity, $\rho$, is
\begin{align}
	\rho(t) &= \rho_{th}(t) + \rho_{ext}(t),
\end{align}
where $\rho_{ext}$ is the external reactivity, and $\rho_{th}$ is the thermal reactivity described by
\begin{align}\label{eq:nuclear:fission:thermal:reactivity}
	\dot \rho_{th}(t) &= -\kappa \dot T_r(t).
\end{align}
Here, $\kappa$ is a proportionality constant, and the reactor temperature, $T_r$, and the temperature in the heat exchanger, $T_{hx}$, are described by the energy balances
\begin{subequations}\label{eq:nuclear:fission:temperatures}
	\begin{align}
		\label{eq:nuclear:fission:temperatures:reactor}
		\dot T_r(t)
		&= \frac{f_r(t)}{m_r} \big(T_{r, in}(t) - T_r(t)\big) + \frac{Q_g(t)}{m_r c_P}, \\
		\label{eq:nuclear:fission:temperatures:heat:exchanger}
		\dot T_{hx}(t)
		&= \frac{f_{hx}(t)}{m_{hx}} \big(T_{hx, in}(t) - T_{hx}(t)\big) \nonumber \\
		&- \frac{k_{hx}}{m_{hx} c_P} \big(T_{hx}(t) - T_c\big).
	\end{align}
\end{subequations}
The mass flow rates through the reactor and heat exchanger,
\begin{align}
	f_r(t) &= f_{hx}(t) = F(t) \rho_s,
\end{align}
are identical. Here, $\rho_s$ is the density of the molten salt, $m_r$ and $m_{hx}$ are the (constant) masses in the reactor and the heat exchanger, $c_P$ is the specific heat capacity at constant pressure, and $k_{hx}$ and $T_c$ are the conductivity and the (constant) temperature of the coolant in the heat exchanger. The thermal energy generated by the fission events is
\begin{align}
	Q_g(t) &= Q_{g, 0} \frac{C_n(t)}{C_{n, 0}},
\end{align}
where $C_{n, 0}$ and $Q_{g, 0}$ are the nominal neutron concentration and the thermal energy generation, respectively. Finally, the inlet temperatures of the reactor and the heat exchanger are
\begin{subequations}
	\begin{align}
		T_{r, in}(t) &= \int_{-\infty}^t \alpha_h(t - s) T_{hx}(s) \incr s, \\
		T_{hx, in}(t) &= \int_{-\infty}^t \alpha_h(t - s) T_r(s) \incr s,
	\end{align}
\end{subequations}
where the kernel $\alpha_h$ is given by~\eqref{eq:hagen:poiseuille:kernel} for a length of $L/2$ and a pressure difference of $\Delta P/2$ because the heat exchanger is located in the middle of the external circulation loop. The manipulated inputs are $\rho_{ext}$ and $\Delta P$.
	\section{Numerical example}\label{sec:numerical:example}
\begin{figure*}[tbh]
	\centering
	\includegraphics[width=\textwidth]{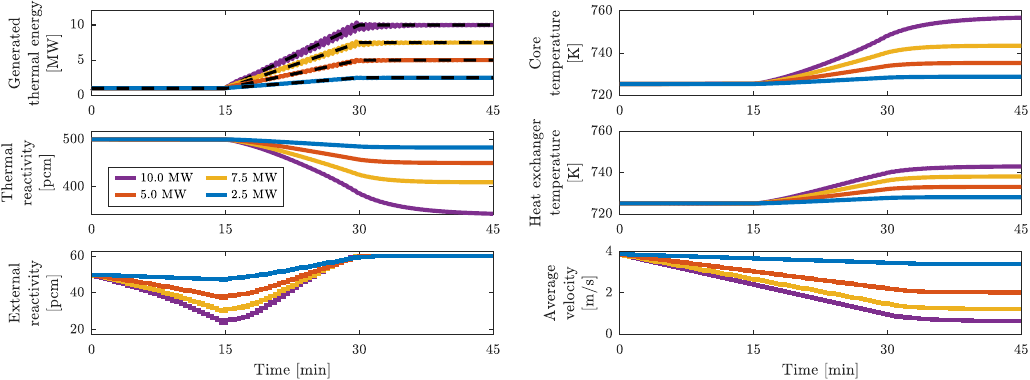}
	\caption{Optimal ramping of the thermal energy generation based on four different time-varying setpoints. Top row: The generated thermal energy and the temperature in the reactor core. Middle row: The thermal reactivity and the temperature in the heat exchanger. Bottom row: The optimal external reactivity and the average velocity corresponding to the optimal pressure difference.}
	\label{fig:optimal:control:comparison}
\end{figure*}

\begin{figure}[tbh]
	\centering
	\includegraphics[width=0.986\linewidth]{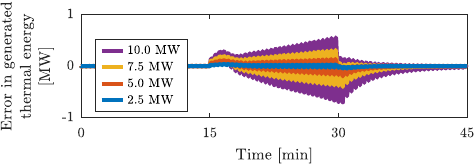}
	\caption{The difference between the generated thermal energy obtained with 1)~the DDEs where the velocity profile is approximated as described in the Appendix and 2)~the approximate system~\eqref{eq:ide} where the delayed variable is linearized.}
	\label{fig:optimal:control:error}
\end{figure}

\begin{figure}[tbh]
	\centering
	\includegraphics[width=0.986\linewidth]{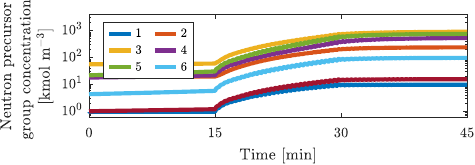}
	\caption{The concentrations of the neutron precursor groups for the simulation shown in Fig.~\ref{fig:optimal:control:comparison} where the thermal energy generation is increased to $10$~MW.}
	\label{fig:optimal:control:neutron:precursor:groups}
\end{figure}

We test the simultaneous approach described in Section~\ref{sec:simultaneous:approach} using the model described in Section~\ref{sec:nuclear:fission} with the parameter values listed in Table~\ref{tab:parameters}. In order to simulate the model (i.e., the true system), we approximate the velocity profile as described in the Appendix where $K = 30$. The resulting system is a set of DDEs with \emph{absolute} time delays, and we use Matlab's \texttt{ddesd} to simulate them.
The control objective is to track a time-varying setpoint that gradually increases the thermal energy generation from 1~MW to 2.5~MW, 5~MW, 7.5~MW, and 10~MW. We use $Q_{g, 0} = 1$~MW and $C_{n, 0} = 1$~kmol~m$^{-3}$, and the initial state is a steady state corresponding to $Q_g = 1$~MW. The corresponding NLPs are solved in 85.6~s, 115.9~s, 97.3~s, and 203.1~s, respectively. The initial guesses of the manipulated inputs are the same as the reference inputs, $u_{-1}$. They are $\rho_{ext} = 50$~pcm and a pressure difference, $\Delta P$, corresponding to an average velocity of 4~m/s. For Hagen-Poiseuille flow, the average velocity is half the maximum velocity. At each point in time, the initial guess of the states is the steady state corresponding to the thermal energy generation. We use control intervals that are $\Delta t = 30$~s long and $M = 1$ time step per control interval. The weight matrix $W_k$ is diagonal with the elements $10^{-2}$~s~pcm$^{-2}$ and $10^2$~s~Pa$^{-1}$.
The results are shown in Fig.~\ref{fig:optimal:control:comparison}, and the setpoints for the thermal energy generation are successfully tracked. The average velocity decreases as the energy generation increases. Consequently, the time delays increase, and the linearization of the delayed variables in~\eqref{eq:delay:approximation} becomes less accurate. This can be seen from Fig.~\ref{fig:optimal:control:error}, where the error in the generated thermal energy is higher for larger increases in the setpoints. Finally, Fig.~\ref{fig:optimal:control:neutron:precursor:groups} shows the concentrations of the neutron precursor groups for the simulation where the setpoint is increased to 10~MW.

\begin{table}
	\centering
	\caption{Model parameter values.}
	\label{tab:parameters}
	\renewcommand{\arraystretch}{1.2}
	\begin{tabular*}{\linewidth}{@{\extracolsep{\fill}} cccccc}
		\hline
		$\lambda_1$~[$s^{-1}$] & $\lambda_2$~[$s^{-1}$] & $\lambda_3$~[$s^{-1}$] & $\lambda_4$~[$s^{-1}$] & $\lambda_5$~[$s^{-1}$] & $\lambda_6$~[$s^{-1}$] \\
		0.0124 & 0.0305 & 0.1110 & 0.3010 & 1.1300 & 3.0000 \\
		\hline
		$\beta_1$~[$-$] & $\beta_2$~[$-$] & $\beta_3$~[$-$] & $\beta_4$~[$-$] & $\beta_5$~[$-$] & $\beta_6$~[$-$] \\
		0.00021 & 0.00141 & 0.00127 & 0.00255 & 0.00074 & 0.00027 \\
		\hline
	\end{tabular*}
	\begin{tabular*}{\linewidth}{@{\extracolsep{\fill}} cccc}
		$\beta$~[$-$] & $\Lambda$ [s] & $c_P$~[MJ~kg$^{-1}$~K$^{-1}$] & $k_{hx}$~[MW~K$^{-1}$]  \\
		0.0065 & $5\cdot 10^{-5}$ & 2$\cdot 10^{-3}$ & 0.5 \\
		\hline
		$\kappa$~[K$^{-1}$] & $\rho_s$~[kg~m$^{-3}$] & $m_r$~[kg] & $m_{hx}$~[kg] \\
		$5\cdot 10^{-5}$ & 2,000 & 10,000 & 2,500 \\
		\hline
		$V$~[m$^3$] & $R$~[m] & $L$~[m] & $T_c$~[K] \\
		0.5 & 0.3 & 30 & 723.15 \\
		\hline
	\end{tabular*}
\end{table}
	\section{Conclusions}\label{sec:conclusions}
In this paper, we present a simultaneous approach for numerical optimal control of DDEs with distributed time delays. We linearize the delayed variables around the current time. The resulting approximate system is a set of implicit differential equations, which we discretize using Euler's implicit method. We use Matlab's fmincon to solve the resulting NLP. Furthermore, we derive the kernel for distributed time delays arising due to flow in a pipe with a nonuniform velocity profile, e.g., Hagen-Poseuille flow. Finally, we demonstrate the efficacy of the simultaneous approach with a setpoint tracking problem involving a molten salt nuclear fission reactor where the molten salt is circulated through a heat exchanger outside of the reactor core.

	\appendix
	We approximate the integrals in the expressions for the volumetric flow rate $F$ in~\eqref{eq:F} and the molar flow rate $f_{i, in}$ in~\eqref{eq:f:in} using a trapezoidal quadrature rule. In this case, the travel time is a function of the radius, i.e., $\tau = \tau(r)$, and the quadrature points are $r_j = j \Delta r$ for $j = 0, \ldots, K$ where $\Delta r = R/K$. The approximate volumetric flow rate is
\begin{align}
	F &\approx 2 \pi \sum_{j=0}^K w_j v(r_j) r_j \Delta r,
\end{align}
and the approximate molar flow rate is
\begin{align}
	f_{i, in}(t) &\approx 2 \pi \sum_{j=0}^K w_j C_i(t - \tau(r_j)) v(r_j) r_j \Delta r,
\end{align}
where $\tau(r_j) = L/v(r_j)$, and the weights are
\begin{align}
	w_j &=
	\begin{cases}
		1, & j = 1, \ldots, K-1, \\
		\frac{1}{2}, & \text{otherwise}.
	\end{cases}
\end{align}
The approximation is similar for the energy balances presented in Section~\ref{sec:nuclear:fission}.

	\bibliographystyle{IEEEtran}
	\bibliography{./ref/References}
\end{document}